\documentclass[12pt]{article}

\def\demo{\noindent{\bf Proof .}}

\newtheorem{proposition}{Proposition}
\newtheorem{definition}{Definition}

\newtheorem{remark}{Remark}

\newtheorem{corollary}{Corollary}

\begin{document}
\begin{center}
{\LARGE\bf \textsc{Factorials as sums}}
\end{center}
\vskip.5truecm
\begin{center}
\vskip 20pt
{\bf Roberto Anglani}\\
{\small\it Dipartimento Interateneo di 
Fisica ``Michelangelo Merlin'', Universit\`{a} degli Studi di Bari,\linebreak Via
 Amendola 173, 70126 Bari,
Italy}\\
{\tt roberto.anglani@ba.infn.it}\\ 
\vskip 10pt
{\bf Margherita Barile\footnote{Partially supported
 by the Italian Ministry of Education, University and 
Research.}}\\
{\small\it Dipartimento di Matematica, Universit\`{a} degli Studi 
di Bari, Via E. Orabona 4, 70125 Bari, Italy}\\
{\tt barile@dm.uniba.it}\\ 
\end{center}
\vskip.5truecm

\section*{Introduction}
In this paper we give an additive representation of the factorial, which can be proven by a simple quick analytical argument.\newline
We also present some generalizations, which are linked, on the one hand to an arithmetical theorem proven by Euler (decomposition of primes as the sum of two squares), and, on the other hand, to modern combinatorics (Stirling numbers). 
\section{From calculus to combinatorics}

\begin{proposition}\label{p1}  For every positive integer $n$, 
$$n!=\sum_{i=0}^n(-1)^{n-i}{n\choose i} i^n.$$
\end{proposition}

\demo For all $x\in{{\bf R}}$ we have the following first order Taylor expansion with Lagrange remainder of the exponential function about 0:

$$e^x=1+x+x^2g(x),$$
\noindent
where $g\in C^{\infty}({\bf R})$. We deduce that, for all $x\in{\bf R}$,

\begin{equation}\label{0}(e^x-1)^n=(x+x^2g(x))^n=x^n+x^{n+1}h(x),\end{equation}
for some $h\in C^{\infty}({\bf R})$. The L.H.S. of (\ref{0}) is 
$$\sum_{i=0}^n(-1)^{n-i}{n\choose i}e^{ix}.$$
\noindent
Its $n$-th derivative is
$$\sum_{i=0}^{n}(-1)^{n-i}{n\choose i}i^ne^{ix},$$
\noindent
whose value at 0 is 
$$\sum_{i=0}^{n}(-1)^{n-i}{n\choose i}i^n.$$
This is equal to the value of the $n$-th derivative of the R.H.S. of (\ref{0}) at 0, which is $n!$. 
This completes the proof.
\vskip\baselineskip\noindent
Proposition \ref{p1} admits an alternative combinatorial proof. Fix  an alphabet of $n$ letters. The number $n!$ counts the ways of forming a word of length $n$ using all $n$ letters (i.e., a word having $n$ different letters), whereas $n^n$ is the total number of words of length $n$ (i.e., of all words, including the possibility of repeated letters).  Hence $n^n-n!$ is the number of words of length $n$ containing at least one repeated letter. Thus we can restate Proposition \ref{p1} in the following equivalent form:
\begin{proposition} The number of words of length $n$ containing at least one repeated letter from an alphabet of $n$ letters is 
\begin{equation}\label{3}-\sum_{i=1}^{n-1}(-1)^{n-i}{n\choose i}i^n.\end{equation}
\end{proposition}
\demo
Recall that, for all $i=1,\dots, n-1$, $i^n$ is the number of words of length $n$ which can be formed from a subset of $i$ letters; ${n\choose i}$ is the number of ways of choosing $i$ letters from our alphabet. Hence the product ${n\choose i}i^n$ counts the words of length $n$ containing at most $i$ different letters, but each word is counted more than once. Precisely, in the sum (\ref{3}), for every integer $k$, $1\leq k\leq n-1$, each word containing exactly $k$ letters is counted one time for each $i$-subset containing these $k$ letters. The number of these subsets is equal to the number of ways in which we can pick $i-k$ letters from the remaining $n-k$ letters of the alphabet, i.e., it is equal to ${n-k\choose i-k}$. Therefore, in (\ref{3}), each word with exactly $k$ letters is counted $-\sum_{i=k}^{n-1}(-1)^{n-i}{n-k\choose i-k}$ times in total. 
Hence the claim is proven once we have shown that
\begin{equation}\label{4b}-\sum_{i=k}^{n-1}(-1)^{n-i}{n-k\choose i-k}=1.\end{equation}
\noindent
Setting $j=i-k$ and making a change of indices, the L.H.S. of (\ref{4b}) becomes
\begin{eqnarray*}-\sum_{j=0}^{n-k-1}(-1)^{n-k-j}{n-k\choose j}&=&-\left(\displaystyle\sum_{j=0}^{n-k}(-1)^{n-k-j}{n-k\choose j}-1\right)\\
&=&-\left((1-1)^{n-k}-1\right)=1,\end{eqnarray*}
\noindent
as was to be shown.
\par\smallskip\noindent
We also have the following different version of Proposition \ref{p1}:
\begin{proposition}\label{p2} For every positive integer $n$,
$$n!=\sum_{i=0}^n(-1)^{n-i}{n\choose i}(i+1)^n.$$
\end{proposition}
\demo
By Proposition \ref{p1} we have
\begin{eqnarray*} (n+1)!&=&\sum_{i=0}^{n+1}(-1)^{n+1-i}{n+1\choose i}i^{n+1}\\
&=&\sum_{i=1}^{n+1}(-1)^{n+1-i}\frac{(n+1)!}{(n+1-i)!i!}i^{n+1}\\
&=&\sum_{i=1}^{n+1}(-1)^{n-(i-1)}\frac{(n+1)n!}{(n-(i-1))!(i-1)!}i^{n}\\
&=&(n+1)\sum_{i=0}^n(-1)^{n-i}\frac{n!}{(n-i)!i!}(i+1)^n\\
&=&(n+1)\left[\sum_{i=0}^n(-1)^{n-i}{n\choose i}i^n+\sum_{i=0}^n(-1)^{n-i}{n\choose i}((i+1)^n-i^n)\right]\\
&=&(n+1)\left\lbrack n!+\sum_{i=0}^n(-1)^{n-i}{n\choose i}((i+1)^n-i^n)\right\rbrack
\end{eqnarray*}
By comparing the first and the last expression in the above sequence of equalities, it follows that
$$\label{2}\sum_{i=0}^n(-1)^{n-i}{n\choose i}((i+1)^n-i^n)=0.$$
\noindent
In view of Proposition \ref{p1} this proves the claim.

\section{Euler's approach}
In the preceding section, Proposition \ref{p2} was obtained from Proposition 1 by replacing $i^n$ with $(i+1)^n$.  We are going to show that, more generally, for every positive integer $n$, the identity 

\begin{equation}\label{4}n!=\sum_{i=0}^n(-1)^{n-i}{n\choose i} (i+k)^n\end{equation}
\noindent
holds for  all nonnegative integers $k$. This will be a consequence of Proposition \ref{p5} below. 
Identity (\ref{4}) was proven by Euler \cite{E1} in a purely arithmetical context. In this section we essentially follow his line of thoughts. 

\begin{definition} Given a sequence $(a_k)_{k\geq 0}$ of real numbers,  its {\it first difference} is defined as the sequence whose $k$-th term is 
$$\Delta^1a_k=a_{k+1}-a_k.$$
\noindent
Recursively, for every integer $n>1$, the {\it n-th difference} is defined as the first difference of $(\Delta^{n-1}a_k)_{k\geq 0}$.
\end{definition}

\begin{proposition}\label{p3} For every positive integer $n$ and all nonnegative integers $k$,
$$\Delta^na_k=\sum_{i=0}^n(-1)^{n-i}{n\choose i}a_{k+i}.$$
\end{proposition}

\demo  Let $k$ be a nonnegative integer. We prove the claim for $k$ by induction on $n$. The claim is obviously true for $n=1$. Now let $n\geq 2$ and suppose the claim true for $n-1$. We prove it for $n$ and $k$. We have
\begin{eqnarray*}\Delta^na_k&=&\Delta^{n-1}a_{k+1}-\Delta^{n-1}a_{k}\\
&=&\sum_{i=0}^{n-1}(-1)^{n-1-i}{n-1\choose i}a_{k+1+i}-\sum_{i=0}^{n-1}(-1)^{n-1-i}{n-1\choose i}a_{k+i}\\
&=&\sum_{i=1}^{n}(-1)^{n-i}{n-1\choose i-1}a_{k+i}-\sum_{i=0}^{n-1}(-1)^{n-1-i}{n-1\choose i}a_{k+i}\\
&=&-(-1)^{n-1}a_k+\sum_{i=1}^{n-1}(-1)^{n-i}\left[{n-1\choose i-1}+{n-1\choose i}\right] a_{k+i}+a_{k+n}\\
&=&(-1)^n{n\choose 0}a_k+\sum_{i=1}^{n-1}(-1)^{n-i}{n\choose i}a_{k+i}+{n\choose n}a_{k+n}\\
&=&\sum_{i=0}^{n}(-1)^{n-i}{n\choose i}a_{k+i},
\end{eqnarray*}
as was to be shown.
\vskip\baselineskip\noindent
Now suppose that for all $k$,
$$a_k=P(k),$$
\noindent
where $P(x)$ is a polynomial  of degree $n$. Let 
$$P(x)=\alpha x^n+\beta x^{n-1}+\mbox{terms of lower degree}\qquad(\alpha,\beta\in{\bf R}, \alpha\ne0).$$ 
\noindent
 Then, for all $k$,
$$\Delta^1a_k=P(k+1)-P(k).$$
\noindent
Note that 
$$P(x+1)-P(x)=\alpha(x+1)^n-\alpha x^n+\beta(x+1)^{n-1}-\beta x^{n-1}+\mbox{ terms of smaller degree }$$
\noindent is  a polynomial of degree $n-1$, its term of highest degree being $n\alpha x^{n-1}$. By finite descending induction it follows that 
$\Delta^na_k$ is a polynomial of degree 0, i.e., a constant polynomial. We have just established 
\begin{proposition}\label{p4} The $n$-th difference of a polynomial sequence of degree $n$ is a constant sequence.
\end{proposition}
For all positive integers $n$, let $_n\Delta^n_k$ denote the $k$-th term of the $n$-th difference of the sequence $(k^n)_{k\geq0}$. According to Proposition \ref{p3}, for all nonnegative integers $k$  we have
\begin{equation}\label{5} _n\Delta^n_k=\sum_{i=0}^{n}(-1)^{n-i}{n\choose i}(i+k)^n.\end{equation}
\noindent
Euler proved
\begin{proposition}\label{p5} Let $n$ be a positive integer. Then for all nonnegative integers $k$, 
$$_n\Delta^n_k=n!.$$
\end{proposition}
\demo  By virtue of Proposition \ref{p4}, the sequence $(_n\Delta^n_k)_{k\geq0}$ is constant. Hence, for all positive integers $n$ and $k$,
$$_n\Delta^n_k={}_n\Delta^n_0 = \sum_{i=0}^{n}(-1)^{n-i}{n\choose i}i^n=n!,$$
\noindent
where the last two equalities follow by (\ref{5}) and Proposition \ref{p1} respectively. This completes the proof.  
\vskip\baselineskip\noindent
From (\ref{5}) and Proposition \ref{p5} we deduce the property that we announced at the beginning of this section. In fact we have the following generalization of Propositions \ref{p1} and  \ref{p2}:
\begin{corollary}\label{coro} Let $n$ be a positive integer. Then for all nonnegative integers $k$,
$$n!=\sum_{i=0}^{n}(-1)^{n-i}{n\choose i}(i+k)^n.$$
\end{corollary}
\begin{remark}{\rm Actually the claim of Corollary \ref{coro} holds for any real number $k$. This can be proven by essentially the same arguments.}\end{remark}
\vskip\baselineskip\noindent
Next we show how Euler applied Proposition \ref{p5} to prove  a famous statement going back to Pierre de Fermat:
\begin{proposition}\label{p6} Every positive prime $p=4n+1$, where $n$ is a positive integer, is the sum of two  squares. 
\end{proposition}
\demo The numbers $k^{2n}$, where $k$ is a positive integer, are not all congruent modulo $p$. Otherwise the terms of the first difference of the sequence  $(k^{2n})_{k\geq0}$ would all be divisible by $p$, and the same would be true for the differences of higher order. But this is impossible: by Proposition \ref{p5}, for all $k$,  
$$_{2n}\Delta^{2n}_k=(2n)!,$$
\noindent
and $p=4n+1$ does not divide any of the factors of $(2n)!$, since these are all smaller than $p$. Hence there are two positive integers $k$ and $h$ such that 
\begin{equation}\label{divide}p\hbox{ does not divide }h^{2n}-k^{2n}.\end{equation}
Evidently we may assume that $1\leq h<k<p$. 
Then $p$ does not divide $h$ and $k$. Consequently, by Fermat's Little Theorem, $h^{p-1}\equiv k^{p-1}\ (\hbox{mod }p)$, so that 
$$p\hbox{ divides }h^{p-1}-k^{p-1}=h^{4n}-k^{4n}=(h^{2n}+k^{2n})(h^{2n}-k^{2n}).$$
\noindent 
Therefore, taking (\ref{divide}) into account, 
$$p\hbox{ divides }h^{2n}+k^{2n}.$$
\noindent According to a former result by Euler (see Propositio 4 in \cite{E2}), if a number is the sum of two  squares, then the same is true for all its prime divisors that do not divide the two squares. Hence $p$ is the sum of two  squares.

\section{Stirling numbers}

Corollary \ref{coro} is an arithmetical generalization of Proposition \ref{p1}. We now give another generalization that is relevant from a combinatorial point of view: it consists in replacing the exponent $n$  with an arbitrary nonnegative integer $\ell$. Following the arguments used in the proof of Proposition \ref{p1}, we can conclude that 
$$\sum_{i=0}^n(-1)^{n-i}{n\choose i} i^{\ell}$$
\noindent
is equal to the $\ell$-th derivative of $(e^x-1)^n$ at 0. This is equal to $\ell!$ times the coefficient $C(\ell, n)$ of $x^{\ell}$ in the Taylor series expansion of $(e^x-1)^n$ at 0.  Since
$$e^x-1=\sum_{i=1}^{\ell}\frac{x^i}{i!}+x^{\ell+1}g(x)$$
\noindent for some $g\in C^{\infty}({\bf R})$, we can compute $C(\ell, n)$ as the coefficient of $x^{\ell}$ in  the following polynomial:
\begin{equation}\label{power}\left(\sum_{i=1}^{\ell} \frac{x^i}{i!}\right)^n.\end{equation}
\noindent
The term of least degree is $x^n$. Therefore, $C(\ell, n)=0$ if $\ell<n$. Suppose that $\ell \geq n$. Then  the term of degree $\ell$ is the sum of all products
\begin{equation}\label{products} \frac{x^{i_1}}{i_1!}\frac{x^{i_2}}{i_2!}\cdots\frac{x^{i_n}}{i_n!}\qquad\mbox{with }i_1+i_2+\cdots+i_n=\ell\end{equation}
\noindent
that can be formed when applying the distributive law to (\ref{power}). For every set $\{i_1,i_2,\dots, i_n\}$ there is one product (\ref{products}) for each possible arrangement of the sequence $i_1,i_2,\dots, i_n$. Call $N(i_1,i_2,\dots, i_n)$ the number of these arrangements. Then it holds:
$$C(\ell, n)=\sum_{\scriptstyle 1\leq i_1\leq i_2\leq\cdots \leq i_n\atop i_1+i_2+\cdots+i_n=\ell}\frac{N(i_1,i_2,\dots, i_n)}{i_1!i_2!\cdots i_n!},$$
\noindent
We have thus proven the following:
\begin{proposition}\label{stirling}
$$\sum_{i=0}^n(-1)^{n-i}{n\choose i} i^{\ell}=
\displaystyle\cases{0 & if\ $0\leq\ell<n$\cr&\cr
\displaystyle\sum_{1\leq i_1\leq i_2\leq\cdots \leq i_n\atop i_1+i_2+\cdots+i_n=\ell}\displaystyle\frac{\ell!N(i_1,i_2,\dots, i_n)}{i_1!i_2!\cdots i_n!}& if\  $\ell\geq n$}$$
\end{proposition}
\noindent
Note that $N(i_1,i_2,\dots, i_n)$ can be interpreted as the number of ways in which a word of length $\ell$ can be subdivided into a sentence with $n$ words of lengths $i_1$, $i_2$, ..., $i_n$. Hence $\ell!N(i_1,i_2,\dots, i_n)$ is the number of sentences of length $\ell$ containing all $\ell$  letters of a given alphabet and formed by $n$ words of lengths $i_1$, $i_2$, ..., $i_n$. If the words are replaced by subsets (i.e., the arrangement of letters in each word does not matter), the number of sentences reduces to 

$$\frac{\ell!N(i_1,i_2,\dots, i_n)}{i_1!i_2!\cdots i_n!}.$$
\noindent
If, in addition, we also assume that the arrangement of the $n$ words does not matter, the number of sentences becomes

$$\frac1{n!}\frac{\ell!N(i_1,i_2,\dots, i_n)}{i_1!i_2!\cdots i_n!}.$$
\noindent
This is the number of ways in which $\ell$ distinguishable elements can be distributed into $n$ non distinguishable subsets with $i_1$, $i_2$, $\dots$, $i_n$ elements respectively.   According to  Proposition \ref{stirling} we have, for all $\ell\geq n$:

$$\sum_{i=0}^n(-1)^{n-i}\frac1{(n-i)!i!} i^{\ell}=\sum_{\scriptstyle 1\leq i_1\leq i_2\leq\cdots \leq i_n\atop i_1+i_2+\cdots+i_n=\ell}\frac1{n!}\frac{\ell!N(i_1,i_2,\dots, i_n)}{i_1!i_2!\cdots i_n!},$$

\noindent and this number is equal to the total number of ways in which $\ell$ elements can be distributed into $n$ non empty subsets. It is well-known  in combinatorics.
\begin{definition}\label{def} {\rm The numbers
 $$S(\ell, n)=\sum_{i=0}^n(-1)^{n-i}\frac1{(n-i)!i!} i^{\ell}\qquad (0<n\leq\ell)$$ are called {\it Stirling numbers (of the second kind)}.} 
\end{definition}
\noindent More details on Stirling numbers can be found in \cite{LW}, Chapter 3, or in \cite{E}.
\begin{remark}{\rm  If in the formula of Definition \ref{def} we replace  $i^{\ell}$ with $(i+1)^{\ell}$, then, as one can easily check, we obtain another Stirling number:}
$$S(\ell+1, n+1)=\sum_{i=0}^n(-1)^{n-i}\frac1{(n-i)!i!} (i+1)^{\ell}.$$
\end{remark}

\section {Negative exponents}
We finally present some identities which extend those of the preceding sections when the exponent is $\ell=-1$. The proofs that are based on easy induction arguments will be omitted.  The first result is related to Proposition \ref{stirling}.
\begin{proposition} For all positive integers $n$
$$\sum_{i=1}^n(-1)^{n-i}{n\choose i}\frac1i=(-1)^{n-1}\sum_{i=1}^n\frac1i.$$
\end{proposition}
\noindent
The term on the R.H.S. of the claim is, up to the sign, the sum of the first $n$ terms of the well-known harmonic series.\newline Next we replace $\displaystyle\frac1i$ with $\displaystyle\frac1{i+k}$. We shall obtain, as in Proposition \ref{p1} and in Corollary \ref{coro}, an expression involving factorials. 
\begin{proposition}\label{k1}For all positive integers $n,k$ it holds
$$\sum_{i=0}^n(-1)^{n-i}{n\choose i}\frac1{i+k}=(-1)^n\frac{n!(k-1)!}{(n+k)!}.$$
\end{proposition}
\demo We prove that the claim is true for all positive integers $n$ when $k=1$.  Once again we use calculus.\newline
We consider the real function $f$ defined, for all $x\in{\bf R}$, by $$f(x)=(1-x)^n=\sum_{i=0}^n(-1)^i{n\choose i}x^i.$$
\noindent The indefinite integral of $f$ is 
$$F_c(x)=\sum_{i=0}^n(-1)^i{n\choose i}\frac{x^{i+1}}{i+1}+c\qquad(c\in{\bf R}),$$
\noindent and the function $F$ defined by
$$F(x)= -\frac1{n+1}(1-x)^{n+1}$$ 
\noindent
is an antiderivative of $f(x)$.  By comparing the constant terms of $F_c(x)$ and $F(x)$,  we find that
$$F(x)=F_{-\frac1{n+1}}(x)=\sum_{i=0}^n(-1)^i{n\choose i}\frac{x^{i+1}}{i+1}-\frac1{n+1}.$$
\noindent
Now
$$0=F(1)=\sum_{i=0}^n(-1)^i{n\choose i}\frac1{i+1}-\frac1{n+1},$$
\noindent
which yields the claim for $k=1$, namely
$$\sum_{i=0}^n(-1)^{n-i}{n\choose i}\frac1{i+1}=(-1)^n\frac{n!}{(n+1)!}.$$
\noindent  The proof of the claim for any given $k\geq2$  is an easy induction on $n\geq0$.

\end{document}